\documentclass[a4paper,11pt,reqno]{amsart}

\usepackage{latexsym,dsfont,amssymb,amsmath,amsthm}

\usepackage[all]{xy}

\chardef\bslash=`\\ 

\numberwithin{equation}{section}

\newtheorem{theorem}{Theorem}[section]
\newtheorem{corollary}[theorem]{Corollary}
\newtheorem{lemma}[theorem]{Lemma}
\newtheorem{proposition}[theorem]{Proposition}

\theoremstyle{remark}
\newtheorem{remark}[theorem]{Remark}

\theoremstyle{definition}
\newtheorem{definition}[theorem]{Definition}

\newcommand\bp{\begin{proof}}
\newcommand\ep{\end{proof}}
\newcommand{\thmref}[1]{Theorem~\ref{#1}}

\newcommand{\proref}[1]{Proposition~\ref{#1}}

\newcommand{\corref}[1]{Corollary~\ref{#1}}

\def\inv{^{-1}}
\newcommand\3[1]{{\mathds #1}}

\newcommand\aaa{\mathfrak a}
\newcommand\bb{\mathfrak b}
\newcommand\mm{\mathfrak m}
\newcommand\pp{\mathfrak p}

\newcommand{\Z}{\mathbb Z}

\newcommand{\Q}{\mathbb Q}
\newcommand{\R}{\mathbb R}

\newcommand\A{\mathbb{A}_{\Q}}
\newcommand\af{\mathbb{A}_{\Q,f}}
\newcommand\ak{{\mathbb A}_K}
\newcommand\akf{{\mathbb A}_{K,f}}
\newcommand\jkf{\akf^*}
\newcommand\ohs{{\hat{\OO}^*}}
\newcommand\Zhat{\hat{\mathbb Z}}

\newcommand{\rr}{{\mathcal R}}
\newcommand\LL{{\mathcal L}}
\newcommand\OO{{\mathcal O}}

\newcommand\gal{\mathcal G}
\newcommand\kab{K^{ab}}

\newcommand\bpmatrix{\begin{pmatrix}}
\newcommand\epmatrix{\end{pmatrix}}

\newcommand\mn{{\operatorname{Mat}_n}}

\newcommand\gln{{\operatorname{GL}_n}}

\newcommand\glnq{{\operatorname{GL}^+_n}(\Q)}
\newcommand\glnr{{\operatorname{GL}^+_n}(\R)}

\newcommand\slnz{{\operatorname{SL}_n}(\Z)}

\newcommand{\Tr}{\operatorname{Tr}}
\newcommand{\supp}{\operatorname{supp}}

\newcommand\enu[1]{\smallskip\newline\makebox[5mm][l]{\rm(#1)}}

\begin{document}
\title{On Bost-Connes type systems for number fields}

\author[M. Laca]{Marcelo Laca$^1$}
\address{Department of Mathematics and Statistics, University of
Victoria, P.O. Box 3045, Victoria, British Columbia, V8W 3P4, Canada.}
\email{laca@math.uvic.ca}
\thanks{$^1$) Supported by the Natural Sciences and Engineering Research
Council of Canada.}

\author[N. S. Larsen]{Nadia S. Larsen$^2$}
\address{Department of Mathematics, University of Oslo,
P.O. Box 1053 Blindern, N-0316 Oslo, Norway.}
\email{nadiasl@math.uio.no}
\thanks{$^2$) Supported by the Research Council of Norway.}

\author[S. Neshveyev]{Sergey Neshveyev$^2$}
\address{Department of Mathematics, University of Oslo,
P.O. Box 1053 Blindern, N-0316 Oslo, Norway.}
\email{sergeyn@math.uio.no}

\begin{abstract}
We give a complete description of the phase transition of the
Bost-Connes type systems for number fields recently introduced by
Connes-Marcolli-Ramachandran and Ha-Paugam. We also introduce a
notion of $K$-lattices and discuss an interpretation of these
systems in terms of $1$-dimensional $K$-lattices.
\end{abstract}

\date{October 30, 2007}

\maketitle

\bigskip

\section*{Introduction}
The generalization of the results of Bost and Connes \cite{bos-con}
to general number fields has received significant attention for more
than ten years, but was only recently formulated in detail as an
explicit problem, see \cite[Problem 1.1]{cmn1}. We paraphrase here
this formulation for easy reference:

{\em Given an algebraic number field $K$, construct a C$^*$-dynamical
system $(\mathcal A, \sigma)$ such that
\begin{itemize}
\item[(i)] the partition function of the system is the Dedekind zeta
function of  $K$; \item[(ii)] the quotient of the idele class group
$C_K$ by the connected component $D_K$ of the identity  acts as
symmetries of the system; \item[(iii)] for each inverse temperature
$0<\beta \leq 1$ there is a unique KMS$_\beta$-state; \item[(iv)]
for each $\beta >1$ the action of the symmetry group $C_K/D_K$ on
the extremal KMS$_\beta$-states is free and transitive; \item[(v)]
there is a $K$-subalgebra $\mathcal A_0$ of $\mathcal A$ such that
the values of the extremal KMS$_\infty$-states on elements of
$\mathcal A_0$ are algebraic numbers that generate the maximal
abelian extension $\kab$ of $K$; and \item[(vi)] the Galois action
of $\gal(\kab/K)$ on these values is realized by the action of $C_K
/ D_K$ on the extremal KMS$_\infty$-states via the class field
theory isomorphism $s\colon C_K / D_K \to \gal(\kab/K)$.
\end{itemize}
 }

Systems with properties (i)-(iv)  have been constructed by several
authors for various classes of number fields,  see \cite[Section
1.4]{con-mar} for a discussion of these constructions and an
extensive list of references. However, the last two properties  have
proven quite elusive. This should not come as a surprise, since a
system satisfying (v) and (vi) has the potential to shed light onto
Hilbert's $12$th problem  about the explicit class field theory of
$K$, although this will ultimately depend on the specific
expressions obtained for the extremal KMS$_\infty$-states and the
generators of the subalgebra~$\mathcal A_0$. Since imaginary
quadratic fields are the only fields beyond $\Q$ for which explicit
class field theory is completely understood, it is natural that they
should be the first case to be solved, and indeed, Connes, Marcolli,
and Ramachandran have produced a complete solution of the problem
for these fields, see \cite[Theorem 3.1]{cmn1}. It should be noticed
also that properties (v) and (vi) are intrinsically related so that
the `right' Galois symmetries and the `right' arithmetic subalgebra
must match each other for the system to have genuine class field
theory content. This failed for instance in the system constructed
in \cite{lvf}, where it was natural to include certain cyclotomic
elements in the arithmetic subalgebra $\mathcal A_0$, but the Galois
action on the corresponding values of extremal KMS$_\infty$ states
did not match the symmetry action of the idele class group on these
elements, see~\cite[Theorem 4.4]{lvf}.

In this paper we study a generalization  of the system from
\cite{cmn1} to all algebraic number fields; this construction is
also isomorphic to a particular case, for Shimura data arising from
a number field, of a general construction due to Ha and Paugam
\cite{ha-pa}.

In Section~\ref{sKMS} we  show how to reduce the study of KMS states
of certain restricted groupoid C$^*$-algebras to measures of the
(unrestricted) transformation group $(G,X)$ satisfying a scaling
condition. Similar results are well known when  $G$ acts freely,
see~\cite{ren}; the key result in this section allows a certain
degree of non-freeness and is motivated by our considerations
in~\cite{finitecm}. We remark that for our applications in
Section~\ref{sBC} we could use instead the earlier results
from~\cite{lac} on crossed products by lattice semigroups. However,
our present results can be applied to a wider class of systems, e.g.
those studied in \cite{lvf} for fields of class number bigger than
one.

Under further assumptions, in \proref{phasetransition}, we give a
natural parametrization of the extremal KMS${}_\beta$-states in
terms of a specific subset of the space~$X$. In
Proposition~\ref{ground} we prove a similar result for ground
states.

In Section~\ref{sBC} we discuss the dynamical system $(\mathcal A,
\sigma)$. We initially construct the C$^*$-algebra $\mathcal A$
along the lines of \cite{con-mar,ha-pa}, using the restricted
groupoid obtained from an action of the group of fractional ideals
on the cartesian product of $\gal(\kab/K) $ by the finite adeles,
balanced over the integral ideles. Since we choose to incorporate
the class field theory isomorphism in the construction, the symmetry
group of the system is $\gal(\kab/K)$ itself. We also indicate that
$\mathcal A$ is a semigroup crossed product of the type discussed in
\cite{lac}, and that, for totally imaginary fields of class number
one, the resulting system is isomorphic to the one constructed  in
\cite{lvf} using Hecke algebras.

Using the results of Section~\ref{sKMS} we show in \thmref{mainthm}
that for an arbitrary number field $K$ the system $(\mathcal A,
\sigma)$ satisfies properties (i) through (iv) above. The
description of the symmetry action and classification of the
KMS${}_\beta$-states generalize the corresponding results of
\cite{cmn1} and complete the initial results of \cite[Section
6]{ha-pa}. The argument goes along familiar lines~\cite{lac,nes},
but we include a complete proof to make the paper self-contained.

Finally, we introduce in Section~\ref{sKlat} a notion of
$n$-dimensional $K$-lattices, which generalizes the $n$-dimensional
$\Q$-lattices from \cite{con-mar} and the $1$-dimensional
$K$-lattices for imaginary quadratic fields $K$ from \cite{cmn1}.
After discussing some of their basic properties we show in
\corref{Klatscal} how $1$-dimensional $K$-lattices are related to
the systems of Section~\ref{sBC}. Therefore these systems can be
introduced without using any class field theory data. This may turn
out to be useful in verifying properties (v) and (vi), as was the
case for $\Q$ and imaginary quadratic fields~\cite{con-mar,cmn1}.

\bigskip
\section{KMS states and measures}\label{sKMS}

Throughout this section we suppose that $G$ is a countable discrete
group acting on a second countable, locally compact, Hausdorff
topological space $X$ and that $Y$ is a clopen subset of $X$
satisfying $G Y = X$. The C$^*$-algebra $C_0(X)\rtimes_r G$ is the
reduced C$^*$-algebra of the transformation groupoid $G\times X$.
Consider the subgroupoid
$$
G\boxtimes Y=\{(g,x)\mid x\in Y,\ gx\in Y\}
$$
and denote by $C^*_r(G\boxtimes Y)$ its reduced C$^*$-algebra. In
other words, $C^*_r(G\boxtimes Y)=\31_Y(C_0(X)\rtimes_r G) \31_Y$,
where $\31_Y$ is the characteristic function of $Y$.

We endow $C^*_r(G\boxtimes Y)$ with the dynamics $\sigma$ associated
to a given homomorphism $N\colon G \to (0,+\infty)$, so
$$
\sigma_t(f)(g,x)=N(g)^{it}f(g,x)
$$
for $t\in\R$ and $f\in C_c(G\boxtimes Y)\subset C^*_r(G\boxtimes
Y)$. Recall that a KMS state for $\sigma$ at inverse temperature
$\beta\in\R$, or $\sigma$-KMS${}_\beta$-state,  is a
$\sigma$-invariant state~$\varphi$ such that $\varphi(ab) =
\varphi(b\sigma_{i\beta}(a))$ for $a$ and $b$ in a set of
$\sigma$-analytic elements with dense linear span.

Denote by $E$ the usual conditional expectation from $C_0(X)
\rtimes_r G$ onto $C_0(X)$: thus with $u_g$ denoting the element in
the multiplier algebra $M(C_0(X)\rtimes_r G)$ corresponding to $g\in
G$, we have $E(fu_g)=f$ if $g=e$ and $0$ otherwise. Observe that the
image under $E$ of the corner $C_r^*(G\boxtimes Y)$ is $C_0(Y)$. By
restriction to $C_0(Y)$, a state~$\varphi$ on $C_r^*(G\boxtimes Y)$
gives rise to a Radon probability measure $\mu$, and conversely, a
Radon probability measure on $Y$ can be extended via the conditional
expectation to a state on $C_r^*(G\boxtimes Y)$. Clearly $(\mu_*
\circ E)|_{C_0(Y)} = \mu_*$, but in general it is not true that
$(\varphi |_{C_0(Y)}) \circ E$ will be the same as~$\varphi$. We
will show that, under certain combined assumptions on $N$ and the
action of $G$ on $X$, the $\sigma$-KMS$_\beta$-states do arise from
their restrictions to $C_0(Y)$, and are thus in one-to-one
correspondence with a class of measures on~$Y$ characterized by a
scaling condition.

\begin{proposition}\label{KMSmeasures}
Under the general assumptions on $G$, $X$, $Y$ and $N$ listed above,
suppose there exist a sequence $\{Y_n\}^\infty_{n=1}$ of Borel
subsets of $Y$ and a sequence $\{g_n\}^\infty_{n=1}$ of elements of
$G$ such that \enu{i} $\cup^\infty_{n=1} Y_n $ contains the set of
points in $Y$ with nontrivial isotropy; \enu{ii} $N(g_n) \neq 1$ for
all $n\ge1$; \enu{iii} $g_n Y_n = Y_n$ for all $n\ge1$.

\smallskip

Then for each $\beta\ne0$ the map $\mu \mapsto \varphi =(\mu_* \circ
E)|_{C_r^*(G\boxtimes Y)}$ is an affine isomorphism between Radon
measures $\mu$ on $X$ satisfying $\mu(Y) = 1$ and the scaling
condition
\begin{equation}\label{normalrescale}
\mu(gZ) = N(g)^{-\beta} \mu(Z)
\end{equation}
for Borel $Z\subset X$,
and  $\sigma$-KMS$_\beta$-states $\varphi$ on $C_r^*(G\boxtimes Y)$.
\end{proposition}

\begin{proof}
It is straightforward to check that any measure satisfying the
scaling condition extends via~$E$ to a KMS$_\beta$-state.

Conversely, let $\varphi$ be a KMS$_\beta$-state. Denote by $\mu$
the probability measure on $Y$ defined by $\varphi|_{C_0(Y)}$.
Applying the KMS-condition to elements of the form
$u_gfu_g^*=f(g^{-1}\,\cdot)$, it is easy to see that
(\ref{normalrescale}) is satisfied for Borel $Z\subset Y$ such that
$gZ\subset Y$.  In particular, $\mu(Y_n)=0$ by conditions (ii) and
(iii).

To show that $\varphi=\mu_*\circ E$ fix $g\ne e$. Let $f\in C_c(Y)$
be such that $g^{-1}(\supp f)\subset Y$. Then $fu_g\in
C_r^*(G\boxtimes Y)$ and we have to prove that $\varphi(fu_g)=0$.

Denote by $Y_g$ the set of points of $Y$ left invariant by $g$. If
$\supp f\cap Y_g=\emptyset$ then we can write $f$ as a finite sum of
functions $h_1h_2$ such that $g(\supp h_1)\cap\supp h_2=\emptyset$.
By the KMS-condition we have
$$
\varphi(h_1h_2u_g)=\varphi(h_2u_gh_1)=\varphi(h_2h_1(g^{-1}\cdot)u_g)=0.
$$
Therefore $\varphi(fu_g)=0$.

Assume now that $\supp f\cap Y_g\ne\emptyset$. As $\mu(Y_n)=0$, by
condition (i) we get $\mu(Y_g)=0$. Hence there exists a norm-bounded
sequence $\{f_n\}_n\subset C_c(Y)$ such that $\supp f_n \cap
Y_g=\emptyset$, $g^{-1}(\supp f_n)\subset Y$ and $f_n\to f$ in
measure $\mu$. Then $\varphi(f_nu_g)=0$. On the other hand, by the
Cauchy-Schwarz inequality,
$$
|\varphi(fu_g)-\varphi(f_n u_g)|\le\varphi(|f-f_n|)^{1/2}
\varphi(u^*_g|f-f_n|u_g)^{1/2}
\le\|f-f_n\|^{1/2}\varphi(|f-f_n|)^{1/2},
$$
whence $\varphi(fu_g)=0$. Therefore $\varphi=\mu_*\circ E$.

To finish the proof it remains to note that the measure $\mu$
extends uniquely to a measure on $X$, which we still denote by
$\mu$, such that (\ref{normalrescale}) is satisfied for all $g\in G$
and Borel $Z\subset X$. Explicitly, we can write
$$
\mu(Z)=\sum_iN(h_i)^\beta\mu(h_iZ\cap Z_i),
$$
where $h_i\in G$ and $Z_i\subset Y$ are such that $X$ is the
disjoint union of the sets $h^{-1}_iZ_i$, see
\cite[Lemma~2.2]{cmgl2}.
 \end{proof}

Our next goal is to classify measures satisfying the scaling
condition. The classification depends on convergence of certain
Dirichlet series. More precisely, when $S$ is a subset of $G$ the
zeta function associated to $S$ is defined to be
\[
\zeta_S(\beta)  := \sum_{s\in S} N(s)^{-\beta} .
\]

\begin{proposition}\label{phasetransition}
Assume the hypotheses of \proref{KMSmeasures}. Let $\beta\ne0$, $S$
be a subset of~$G$, and $Y_0 \subset Y$ a nonempty Borel set such
that \enu{i} $gY_0 \cap Y_0 = \emptyset$ for $g\in G \setminus
\{e\}$; \enu{ii} $SY_0 \subset Y$;\enu{iii} if $gY_0\cap
Y\ne\emptyset$ then $g\in S$; \enu{iv} $Y\setminus S
U\subset\cup_nY_n$ for every open set $U$ containing $Y_0$; \enu{v}
$\zeta_S(\beta)<\infty$.

\smallskip

Then \enu{1}  the map $\varphi=\mu_*\circ
E\mapsto\zeta_S(\beta)\mu|_{Y_0}$ is an affine isomorphism between
the $\sigma$-KMS$_\beta$-states on $C^*_r(G\boxtimes Y)$ and the
Borel probability measures  on $Y_0$; the inverse map is given by
$\nu\mapsto\mu_*\circ E$, where the measure $\mu$ on~$Y$ is defined
by
\begin{equation} \label{eextension}
\mu(Z)=\zeta_S(\beta)^{-1}\sum_{s\in S}N(s)^{-\beta}\nu(s^{-1}Z\cap Y_0);
\end{equation}
\enu{2} if $\mu$ is the measure on $Y$ defined by a probability
measure $\nu$ on $Y_0$ by \eqref{eextension}, and $H_S$ is the
subspace of $L^2(Y, d \mu)$ consisting of functions $f$ such that
$f(sy)=f(y)$ for $y\in Y_0$ and $s\in S$, then for $f\in H_S$ we
have
\begin{equation} \label{enorm}
\|f\|^2_2=\zeta_S(\beta)\int_{Y_0}|f(y)|^2d\mu(y);
\end{equation}
furthermore, the
orthogonal projection $P\colon L^2(Y,d\mu)\to H_S$ is given by
\begin{equation} \label{epro0}
Pf|_{Sy}=\zeta_S(\beta)^{-1}\sum_{s\in S}N(s)^{-\beta}f(sy)\ \
\text{for}\ \ y\in Y_0.
\end{equation}
\end{proposition}

\bp By \proref{KMSmeasures} any KMS$_\beta$-state is determined by a
Radon measure $\mu$ such that $\mu(Y)=1$ and $\mu$ satisfies the
scaling condition \eqref{normalrescale}. By assumptions (i) and
(ii), for such a measure $\mu$ we have
$$
1\ge\mu(SY_0)=\sum_{s\in
S}N(s)^{-\beta}\mu(Y_0)=\zeta_S(\beta)\mu(Y_0).
$$
On the other hand, as $\mu(Y_n)=0$, by assumption (iv) we have
$$
1=\mu(Y)\le\mu(SU)\le\sum_{s\in S}\mu(sU)=\zeta_S(\beta)\mu(U)
$$
for any open set $U$ containing $Y_0$. By regularity of the measure
we conclude that $\zeta_S(\beta)\mu(Y_0)\ge1$, and hence
$\zeta_S(\beta)\mu(Y_0)=1$. It follows that $SY_0$ is a subset of
$Y$ of full measure. Since $\mu$ satisfies the scaling condition, we
conclude that $\mu$ is completely determined by its restriction to
$Y_0$.

To finish the proof of (1) we have to construct the inverse map. Let
$\nu$ be a Borel measure on $Y_0$ with $\nu(Y_0)=1$. Similarly to
the proof of Proposition~\ref{KMSmeasures} define a measure $\mu$ on
$X$ by
$$
\mu(Z)=\zeta_S(\beta)^{-1}\sum_{g\in G}N(g)^\beta\nu(gZ\cap Y_0)\ \
\hbox{for Borel}\ \ Z\subset X.
$$
Then $\zeta_S(\beta)\mu$ extends $\nu$ by assumption (i) and
satisfies \eqref{normalrescale}. Furthermore, by assumptions (ii)
and (iii) we have $gY\cap Y_0\ne\emptyset$ if and only if $g^{-1}\in
S$, and in the latter case $Y_0\subset gY$. It follows that for
$Z\subset Y$ we have \eqref{eextension}. In particular,
$\mu(Y)=\zeta_S(\beta)^{-1}\sum_{s\in S}N(s)^{-\beta}\nu(Y_0)=1$.

\smallskip

Turning to the proof of (2), suppose $\mu$ is the measure on $Y$
defined by a probability measure $\nu$ on $Y_0$ by
\eqref{eextension} and  recall that we have already shown that
$SY_0$ is a subset of $Y$ of full $\mu$-measure. Then (2) is a
particular case of \cite[Lemma~2.9]{cmgl2}. For the reader's
convenience we sketch a proof.

Equality \eqref{enorm} follows from the identity
$$
\int_{sY_0}|f|^2d\mu(y)=N(s)^{-\beta}\int_{Y_0}|f(s\,\cdot)|^2d\mu(y),
$$
valid for $f\in L^2(Y,d\mu)$, on summing over $s\in S$. Furthermore,
as
$$
\sum_{s\in S}N(s)^{-\beta}|f(sy)|^2\ge\zeta_S(\beta)
\left|\zeta_S(\beta)^{-1}\sum_{s\in S}N(s)^{-\beta}f(sy)\right|^2
$$
the above identity and \eqref{enorm} show that the operator $T$ on
$L^2(Y,d\mu)$ defined by the right hand side of \eqref{epro0} is a
contraction. Since $Tf=f$ for $f\in H_S$, and the image of~$T$
is~$H_S$, we conclude that $T=P$. \ep

In our applications the set $S$ will be a subsemigroup of $\{g\in
G\mid N(g)\ge1\}$ and $Y_0$ the complement of the union of the sets
$gY$, $g\in S\setminus\{e\}$.

\smallskip

We next give a similar classification of ground states. Recall that
a $\sigma$-invariant state $\varphi$ is called a ground state if the
holomorphic function $z\mapsto\varphi(a\sigma_z(b))$ is bounded on
the upper half-plane for $a$ and $b$ in a set of $\sigma$-analytic
elements spanning a dense subspace. If a state $\varphi$ is a
weak$^*$ limit point of a sequence of states $\{\varphi_n\}_n$ such
that $\varphi_n$ is a $\sigma$-KMS$_{\beta_n}$-state and
$\beta_n\to+\infty$ as $n\to\infty$, then $\varphi$ is a ground
state. Such ground states are called
$\sigma$-KMS$_\infty$-states~\cite{con-mar}.

\begin{proposition}\label{ground}
Under the general assumptions on $G$, $X$, $Y$ and $N$ listed before
Proposition~\ref{KMSmeasures}, define $Y_0=Y\setminus\cup_{\{g\colon
N(g)>1\}}gY$. Assume $Y_0$ has the property that if $gY_0\cap
Y_0\ne\emptyset$ for some $g\in G$ then $g=e$. Then the map
$\mu\mapsto\mu_*\circ E$ is an affine isomorphism between the Borel
probability measures on~$Y$ supported on $Y_0$ and the ground states
on $C^*_r(G\boxtimes Y)$.
\end{proposition}

\bp Assume first that $\mu$ is a probability measure on $Y$
supported on $Y_0$, $\varphi=\mu_*\circ E$. If $a=f_1u_g$ and
$b=f_2u_h$ with $g^{-1}(\supp f_1),h^{-1}(\supp f_2)\subset Y$, then
$E(a\sigma_z(b))$ is nonzero only if $h=g^{-1}$. In the latter case
the function $\varphi(a\sigma_z(b))=N(g)^{-iz}\varphi(ab)$ is
clearly bounded on the upper half-plane if $N(g)\le1$. So assume
$N(g)>1$. As $u_gf_2u_g^{-1}=f_2(g^{-1}\,\cdot)$ is supported on
$gY$, we see that the support of $f_1f_2(g^{-1}\,\cdot)$ is
contained in $Y\setminus Y_0$, whence $\varphi(a\sigma_z(b))=0$.

\smallskip

Conversely, assume $\varphi$ is a ground state. Let $\mu$ be the
probability measure on $Y$ defined by $\varphi|_{C_0(Y)}$. Take an
element $g\in G$ with $N(g)>1$. If $f\in C_c(Y\cap g^{-1}Y)$ is
positive, $a=u_gf^{1/2}$ and $b=f^{1/2}u_{g^{-1}}$, then the
function $z\mapsto\varphi(a\sigma_z(b))$ can be bounded on the upper
half-plane only if it is identically zero. Therefore
$\varphi(f(g^{-1}\,\cdot))=0$. Hence $\mu(gY\cap Y)=0$. Thus $\mu$
is supported on $Y_0$.

It remains to show that $\varphi(fu_g)=0$ for all $g\ne e$ and $f\in
C_c(Y)$ with $g^{-1}(\supp f)\subset Y$. If $x\in\supp f\cap Y_0$
then $g^{-1}x\notin Y_0$ by our assumptions on $Y_0$. Hence there
exists $h\in G$ with $N(h)>1$ such that $g^{-1}x\in hY$. This shows
that the sets $Y\setminus Y_0$ and $ghY$ with $N(h) >1$ form an open
cover of $\supp f$. Using a partition of unit subordinate to this
cover we decompose $f$ into a finite sum of functions with supports
contained in these sets. Therefore we may assume that either $\supp
f\subset Y\setminus Y_0$ or $g^{-1}(\supp f)\subset hY$ for some $h$
with $N(h)>1$. In the first case we have $\varphi(fu_g)=0$ as $\mu$
is supported on $Y_0$. In the second case write $f$ as a
product~$f_1f_2$ of continuous functions with the same support,
letting e.g. $f_1=|f|^{1/2}$ and $f_2=f|f|^{-1/2}$. Consider the
elements $a=f_1u_{gh}$ and $b=f_2(gh\cdot)u_{h^{-1}}$ of
$C^*_r(G\boxtimes Y)$, so that $fu_g=ab$. Since $N(h)>1$, the
function $z\mapsto\varphi(a\sigma_z(b))$ can be bounded on the upper
half-plane only if it is identically zero. Therefore
$\varphi(fu_g)=0$. \ep

\bigskip
\section{Bost-Connes systems for number fields}\label{sBC}

Suppose $K$ is an algebraic number field with subring of integers
$\OO$. Recall some notation. Denote by~$V_K$ the set of places of
$K$, and by $V_{K,f}\subset V_K$ the subset of finite places. For
$v\in V_K$ denote by~$K_v$ the corresponding completion of $K$. If
$v$ is finite, let $\OO_v$ be the closure of $\OO$ in $K_v$. The
ring of finite integral adeles is $\hat\OO=\prod_{v\in
V_{K,f}}\OO_v$, and $\akf=K\otimes_{\OO}\hat\OO$ is the ring of
finite adeles. Denoting by $K_\infty=\prod_{v|\infty}K_v$ the
completion of $K$ at all infinite places, we get the ring
$\ak=K_\infty\times\akf$ of adeles. The idele group is $I_K=\ak^*$.

Consider the topological space $\gal(\kab / K)\times \akf$, where
$\gal(\kab/K)$ is the Galois group of the maximal abelian extension
of $K$. On this space there is an action of the group $\jkf$ of
finite ideles, via the Artin map $s\colon I_K \to \gal(\kab/K)$ on
the first component and via multiplication on the second component:
\[
j(\gamma, m) = (\gamma s(j)\inv , jm)  \ \ \hbox{for}\ \ j\in \jkf,
\ \ \gamma\in \gal(\kab/K), \ \ m\in \akf.
\]
Following~\cite{cmn1} we consider the quotient space
\[
X := \gal(\kab / K)\times_{\ohs} \akf
\]
in which the direct product  is balanced over  the compact open
subgroup of integral ideles $\ohs \subset \jkf$, in the sense that
one takes the quotient by the action given by $u (\gamma, m) =
(\gamma s(u)\inv , um)$ for $u\in \ohs$. This enables a quotient
action of the quotient group $ \jkf/\ohs$, which is isomorphic to
the (discrete) group $J_K$ of fractional ideals in $K$. We remark
that the space $X$ is isomorphic to the one that arises from the
construction of Ha and Paugam when applied to the Shimura data
associated to a number field, see  \cite[Definition 5.5]{ha-pa}.

Finally we restrict to  the clopen subset $Y := \gal(\kab /
K)\times_\ohs \hat\OO$ of $X$, and we consider the dynamical system
$(C^*_r(J_K \boxtimes Y),  \sigma)$, in which the dynamics $\sigma$
is defined in terms of the absolute norm $N\colon J_K \to
(0,+\infty)$. Denote by $J_K^+\subset J_K$ the subsemigroup of
integral ideals, and recall that the norm of such an ideal $\aaa$ is
given by $|\OO/\aaa|$. Remark that by Theorem 2.1 and Theorem 2.4 of
\cite{minautex}  the corner  $C^*_r(J_K\boxtimes Y)= \31_Y ( C_0(X)
\rtimes J_K ) \31_Y$ is the semigroup crossed product $C(Y) \rtimes
J_K^+$.

In this situation the zeta function of the semigroup $J_K^+$ is
precisely the Dedekind zeta function $\zeta_K(\beta) = \sum_{\aaa\in
J_K^+} N(\aaa)^{-\beta}$; it converges for $\beta>1$ and diverges
for $\beta\in(0,1]$.

\begin{theorem} \label{mainthm}
For the system $(C(\gal(\kab / K)\times_\ohs \hat\OO)\rtimes
J_K^+,\sigma)$ we have: \enu{i} for $\beta<0$ there are no
KMS$_\beta$-states; \enu{ii} for each $0< \beta \leq 1$ there is a
unique KMS$_\beta$-state; \enu{iii} for each  $1< \beta<\infty$ the
extremal KMS$_\beta$-states are indexed by $Y_0 :=  \gal(\kab /
K)\times_\ohs \ohs \cong \gal(\kab / K)$, with the state
corresponding to $w\in Y_0$ given by
\begin{equation} \label{eKMSext}
\varphi_{\beta,w} (f) = \frac{1}{\zeta_K(\beta) }\sum_{\aaa\in
J_K^+} N(\aaa)^{-\beta} f(\aaa w)\ \ \hbox{for}\ \ f\in C(\gal(\kab
/ K)\times_\ohs \hat\OO);
\end{equation}
\enu{iv} the ground states of the system are KMS$_\infty$-states,
and the extremal ground states are indexed by~$Y_0$, with the state
corresponding to $w\in Y_0$ given by $\varphi_{\infty,w}(f)=f(w)$.
\end{theorem}

\bp We apply Proposition~\ref{KMSmeasures} to $G=J_K$, $X= \gal(\kab
/ K)\times_{\ohs} \akf$ and $Y=\gal(\kab / K)\times_\ohs\hat\OO$. If
the image of a point $(\alpha,a)\in\gal(\kab/K)\times\akf$ in $X$
has nontrivial isotropy then $a_v=0$ for some~$v$, since this is
true already for the action of $J_K=\jkf/\ohs$ on $\akf/\ohs$.
Therefore for the sequence $\{(g_n,Y_n)\}_n$ we can take the pairs
$(\pp_v,Y_v)$ indexed by the finite places $v$, where $\pp_v$ is the
prime ideal of $\OO$ corresponding to $v$ and $Y_v\subset Y$
consists of the images in $X$ of all pairs
$(\alpha,a)\in\gal(\kab/K)\times\hat\OO$ with $a_v=0$. By
Proposition~\ref{KMSmeasures} we conclude that the
KMS$_\beta$-states for $\beta\ne0$ correspond to the measures $\mu$
on $X$ such that $\mu(Y)=1$ and $\mu$ satisfies the scaling
condition \eqref{normalrescale}.

\smallskip

Clearly there are no such measures for $\beta<0$, since otherwise
the inclusion $\aaa Y\subset Y$ would imply $N(\aaa)^{-\beta}\le1$.
This proves (i).

\smallskip

To prove part (iii) notice that $S=J_K^+$ and $Y_0=\gal(\kab /
K)\times_\ohs \ohs\subset Y$ satisfy conditions (i), (ii) and (iii)
of Proposition~\ref{phasetransition}. In order to verify condition
(iv) let $A\subset V_{K,f}$ be a finite set and denote by $\OO_A$
the product of $\OO_v$ over $v\in A$, and by $\hat\OO_A$ the product
of $\OO_v$ over $v\notin A$, so that $\hat\OO=\OO_A\times\hat\OO_A$.
Consider the open subset
$$
W_A=\gal(\kab /
K)\times_{\ohs}\left(\OO_A^*\times\hat\OO_A\right)
$$
of $Y$. The intersection of these sets over all finite $A$ coincides
with $Y_0$. Since $Y$ is compact and the sets $W_A$ are closed, it
follows that any neighborhood of $Y_0$ contains $W_A$ for some $A$.
The complement of $J_K^+W_A$ in $Y$ consists of  the images of
points $(\alpha,a)\in\gal(\kab/K)\times\hat\OO$ such that $a_v=0$
for some $v\in A$, so it is covered by the sets $Y_v$, $v\in A$,
introduced above. Thus by Proposition~\ref{phasetransition} for each
$\beta>1$ there is a one-to-one affine correspondence between the
KMS$_\beta$-states and the probability measures on~$Y_0$. In particular,
the extremal KMS$_\beta$-states correspond to points of $Y_0$ via
\eqref{eKMSext}, which is a particular case of \eqref{eextension}.
This finishes the proof of part (iii).

\smallskip

Part (iv) follows from (iii) and Proposition~\ref{ground}.

\smallskip

Turning to (ii), we shall first explicitly construct for each
$\beta\in(0,1]$ a measure $\mu_\beta$ on $X$ such that
$\mu_\beta(Y)=1$ and $\mu_\beta$ satisfies the scaling condition
\eqref{normalrescale}. Define $\mu_\beta$ as the push-forward of the
product measure $\mu_\gal\times\prod_{v\in V_{K,f}}\mu_{\beta,v}$ on
$\gal(\kab/K)\times\akf$, where $\mu_\gal$ is the normalized Haar
measure on $\gal(\kab/K)$ and the measures~$\mu_{\beta,v}$ on $K_v$
are defined as follows. The measure $\mu_{1,v}$ is the additive Haar
measure on $K_v$ normalized by $\mu_{1,v}(\OO_v)=1$. The measure
$\mu_{\beta,v}$ is defined so that it is equivalent to $\mu_{1,v}$
and
$$
\frac{d\mu_{\beta,v}}{d\mu_{1,v}}(a)
=\frac{1-N(\pp_v)^{-\beta}}{1-N(\pp_v)^{-1}}\|a\|_v^{\beta-1},
$$
where $\|\cdot\|_v$ is the normalized valuation in the class $v$, so
$\|\pi\|_v=N(\pp_v)^{-1}$ for any uniformizing parameter $\pi
\in\pp_v$. Equivalently, $\mu_{\beta,v}$ is the unique measure on
$K_v$ such that the restriction of $\mu_{\beta,v}$ to $\OO^*_v$ is
the (multiplicative) Haar measure normalized by
$\mu_{\beta,v}(\OO^*_v)=1-N(\pp_v)^{-\beta}$, and $\mu_{\beta,v}(\pi
Z)=N(\pp_v)^{-\beta}\mu_{\beta,v}(Z)$.

To show that the measure $\mu_{\beta}$ is unique it suffices to show
that the action of $J_K$ on $(X,\mu)$ is ergodic for every measure
$\mu$ on $X$ such that $\mu(Y)=1$ and $\mu$ satisfies the scaling
condition \eqref{normalrescale}. Indeed, the set of such measures is
affine, so if all measures are ergodic the set must consist of one
point.

Equivalently, we have to show that the subspace $H$ of $L^2(Y,d\mu)$
of $J_K^+$-invariant functions consists of scalars. Denote by $P$
the projection onto this space. It is enough to compute how $P$ acts
on the pull-backs of functions on $\gal(\kab/K)\times_\ohs\OO_A$ for
finite $A\subset V_{K,f}$. Denote by $J^+_{K,A}$ the unital
subsemigroup of $J^+_K$ generated by $\pp_v$, $v\in A$. Modulo a set
of measure zero $\gal(\kab/K)\times_\ohs\OO_A$ is the union of the
sets $\aaa(\gal(\kab/K)\times_\ohs\OO_A^*)$, $\aaa\in J^+_{K,A}$.
The compact set $\gal(\kab/K)\times_\ohs\OO_A^*$ is a group
isomorphic to $\gal(\kab/K)/s(\ohs_A)$. Therefore it suffices to
compute $Pf$ for the pull-back $f$ of the function
$$
\gal(\kab/K)\times_\ohs\OO_A\ni a\mapsto
\begin{cases}
\tilde\chi(\aaa^{-1}a), &\ \hbox{if}\ \
a\in\aaa(\gal(\kab/K)\times_\ohs\OO_A^*),\\
0, &\ \hbox{otherwise},
\end{cases}
$$
where $\tilde\chi$ is a character of
$\gal(\kab/K)\times_\ohs\OO_A^*$. The character $\tilde\chi$ is
defined by a Dirichlet character $\chi\mod\mm$ with $\mm\in
J^+_{K,A}$.

For a finite set $B\subset V_{K,f}$ denote by $P_B$ the projection
onto the subspace $H_B\subset L^2(Y,d\mu)$ of $J^+_{K,B}$-invariant
functions. Apply Proposition~\ref{phasetransition}(2) with
$G=J_{K,B}:=(J^+_{K,B})^{-1}J^+_{K,B}$, $S=J^+_{K,B}$ and
$Y_0=W_B=\gal(\kab / K)\times_{\ohs}(\OO_B^*\times\hat\OO_B)$. Note
that $\zeta_{J^+_{K,B}}(\beta)=\prod_{v\in
B}(1-N(\pp_v)^{-\beta})^{-1}$. Furthermore, for $\bb\in J^+_{K,B}$
the set $\bb W_B$ intersects the support of $f$ only if $\aaa|\bb$
and the ideals $\aaa$ and $\bb\aaa^{-1}$ are relatively prime, or
equivalently, $\aaa\in J^+_{K,B}$ and $\bb\in\aaa J^+_{K,B\setminus
A}$. Therefore, assuming $A\subset B$, by \eqref{epro0} we get
\begin{align}
P_Bf|_{J^+_{K,B}a}
&=\prod_{v\in B}(1-N(\pp_v)^{-\beta}) \sum_{\mathfrak{c}\, \in
J^+_{K,B\setminus A}}N(\aaa\mathfrak{c})^{-\beta}\tilde\chi(\mathfrak{c} a)\notag \\
&=N(\aaa)^{-\beta}\tilde\chi(a){\prod_{v\in
B}(1-N(\pp_v)^{-\beta})}\sum_{\mathfrak{c} \, \in
J^+_{K,B\setminus A}}N(\mathfrak{c})^{-\beta}\chi(\mathfrak{c})\notag \\
&=N(\aaa)^{-\beta}\tilde\chi(a)\frac{\prod_{v\in
B}(1-N(\pp_v)^{-\beta})}{\prod_{v\in B\setminus
A}(1-\chi(\pp_v)N(\pp_v)^{-\beta})}\notag
\end{align}
for $a\in W_B$. If $\chi$ is trivial we see that $P_Bf$ is constant,
and hence so is $Pf$. On the other hand, for nontrivial~$\chi$ we
get
$$
\|Pf\|_2=\lim_B\|P_Bf\|_2=N(\aaa)^{-\beta}\lim_B\frac{\prod_{v\in
B}|1-N(\pp_v)^{-\beta}|}{\prod_{v\in B\setminus
A}|1-\chi(\pp_v)N(\pp_v)^{-\beta}|}.
$$
The right hand side divided by $N(\aaa)^{-\beta}$ is an increasing
function in $\beta$ on $(0,+\infty)$. For $\beta>1$ it equals
$|L(\chi,\beta)|/\zeta_K(\beta)$. As $L(\chi,\cdot)$ does not have a
pole at $1$, see e.g. \cite[Lemma 13.3]{neu}, we conclude that the
right hand side is zero for $\beta\in(0,1]$. Therefore in either
case we see that $Pf$ is constant. \ep

\begin{remark}
\mbox{\ } \enu{i} There is an obvious action of the Galois group
$\gal(\kab / K)$ of the maximal abelian extension of~$K$ on  $Y$,
given by $\alpha (\gamma, m) = (\alpha\gamma, m)$, and this  gives
rise to an  action of $\gal(\kab / K)$ as symmetries of
$(C^*_r(J_K\boxtimes Y), \sigma)$. This action is clearly free and
transitive on the set $Y_0$ parametrizing the extreme
KMS$_\beta$-states. \enu{ii} It is known~\cite{ha-pa} and easy to
check that the partition function of our system is the Dedekind zeta
function. More precisely, if $\varphi_{\beta,w}$ is an extremal
KMS$_\beta$-state for some $\beta>1$, and $H_{\beta,w}$ is the
generator of the canonical one-parameter unitary group implementing
$\sigma$ in the GNS-representation of $\varphi_{\beta,w}$, then
$\Tr(e^{-\beta H_{\beta,w}})=\zeta_K(\beta)$. \enu{iii} For totally
imaginary fields of class number one the C$^*$-algebra $C^*_r(J_K
\boxtimes Y)$ described above is isomorphic to the Hecke
C$^*$-algebra $C^*(\Gamma_K; \Gamma_\OO)$ studied in \cite{lvf}. To
see this, observe first that $\gal(\kab/K)\cong
\akf^*/\overline{K^*}\cong\ohs/\overline{\OO^*}$. It follows that $Y
= \gal(\kab/K) \times_\ohs \hat\OO $ can be identified
with~$\hat\OO/{\overline{\OO^*}}$. From \cite[Definition 2.2]{lvf}
and the ensueing discussion, multiplication by an extreme inverse
different transforms this identification into a homeomorphism of the
orbit  space $\Omega = \mathcal D\inv /{\overline{\OO^*}}$ and  $Y$.
It is then easy to check that the multiplicative action of $J_K^+
\cong \OO^\times/\OO^*$ on $C(\Omega) $ described  in
\cite[Proposition 2.4]{lvf} corresponds to the action of $J_K^+
\cong \jkf/ \ohs$ inherited by $C(Y)$ from the original
transformation group. By \cite[Theorem 2.5]{lvf} it follows that the
Hecke C$^*$-algebra $C^*(\Gamma_K; \Gamma_\OO)$ is isomorphic to
$C(Y) \rtimes J_K^+ \cong C^*_r(J_K\boxtimes Y)$. The isomorphism
respects the semigroup of isometries and thus the dynamics arising
from the norm, but the Galois group action is changed via the
balancing over $\ohs$, and this resolves the incompatibility pointed
out in \cite[Theorem 4.4]{lvf}.

For higher class numbers the Hecke C$^*$-algebra constructed in
\cite{lvf} is a semigroup crossed product by the semigroup of
principal ideals so it is essentially different from the one studied
here.

\end{remark}

\bigskip
\section{$K$-lattices} \label{sKlat}

In this section we define $n$-dimensional $K$-lattices and interpret
the BC-systems for number fields in terms of $K$-lattices.

Recall the following definition given by Connes and
Marcolli~\cite{con-mar}. An $n$-dimensional $\Q$-lattice is a pair
$(L,\varphi)$, where $L\subset\R^n$ is a lattice and
$\varphi\colon\Q^n/\Z^n\to\Q L/L$ is a homomorphism. The notion of a
$1$-dimensional $K$-lattice for an imaginary quadratic field $K$ is
analyzed in \cite{cmn1}. In what follows we generalize $K$-lattices
to arbitrary number fields and dimensions. We refer to \cite{KCoMa}
for a related discussion of the function fields case, see
also~\cite{J}.

Recall that we denote by $K_\infty$ the completion of $K$ at all
infinite places, so $K_\infty\cong\R^{[K\colon\Q]}$ as a topological
group under addition. By an $n$-dimensional $\OO$-lattice we mean a
lattice $L$ in $K_\infty^n$ such that $\OO L=L$.

\begin{definition}
An $n$-dimensional $K$-lattice is a pair $(L,\varphi)$, where
$L\subset K_\infty^n$ is an $n$-dimensional $\OO$-lattice and
$\varphi\colon K^n/\OO^n\to KL/L$ is an $\OO$-module map.
\end{definition}

The simplest example of an $n$-dimensional $\OO$-lattice is $\OO^n$.
Since $K^n=\Q\OO^n$, any two finitely generated $\OO$-submodules of
$K^n$ of rank $n$ are commensurable, in particular, any such module
is an $\OO$-lattice. Furthermore, a submodule of $K^n$ of rank $m<n$
is an abelian group of rank $m[K\colon\Q]$, so it cannot be a
lattice in $K^n_\infty$. Thus for submodules of $K^n$ we get the
usual definition of an $\OO$-lattice: an $\OO$-submodule $M\subset
K^n$ is an $n$-dimensional $\OO$-lattice if and only if it is
finitely generated and has rank $n$.

We now want to give a parametrization of the set of $n$-dimensional
$\OO$-lattices. For this recall that there exists a one-to-one
correspondence between finitely generated $\OO$-submodules of~$K^n$
of rank $n$ and $\hat\OO$-submodules $\LL=\prod_{v\in
V_{K,f}}L_v\subset\akf^n$ such that $L_v$ is a compact open
$\OO_v$-submodule of $K_v^n$ with $L_v=\OO_v^n$ for all but a finite
number of places $v$. Namely, starting from an $\OO$-lattice
define~$\LL$ as its closure. The inverse map is
$\LL\mapsto\cap_v(L_v\cap K^n)$. Hence, given an element
$s=(s_\infty,s_f)\in\gln(\ak)=\gln(K_\infty)\times\gln(\akf)$, we
get an $\OO$-lattice $s_f\hat\OO^n\cap K^n$ in $K^n$, and then an
$\OO$-lattice $s^{-1}_\infty (s_f\hat\OO^n\cap K^n)$ in
$K^n_\infty$.

\begin{lemma} \label{olat}
The map $\gln(\ak)\ni s\mapsto s^{-1}_\infty (s_f\hat\OO^n\cap K^n)$
induces a bijection between
$$
\gln(K)\backslash\gln(\ak)/\gln(\hat\OO)
$$
and the set of $n$-dimensional $\OO$-lattices.
\end{lemma}

\bp It is easy to see that the map from
$\gln(K)\backslash\gln(\ak)/\gln(\hat\OO)$ to $\OO$-lattices is
well-defined. To see that it is injective, assume $r^{-1}_\infty
(r_f\hat\OO^n\cap K^n)=s^{-1}_\infty (s_f\hat\OO^n\cap K^n)$ for
some $r,s\in\gln(\ak)$. Multiplying by $K$ we get $r^{-1}_\infty
K^n=s^{-1}_\infty K^n$, so $g:=s_\infty r^{-1}_\infty\in\gln(K)$.
Taking the closure we get from $g(r_f\hat\OO^n\cap
K^n)=s_f\hat\OO^n\cap K^n$ that $gr_f\hat\OO^n=s_f\hat\OO^n$. Hence
$gr_fu=s_f$ for some $u\in\gln(\hat\OO)$. Since also
$gr_\infty=s_\infty$, this means that $s$ is in a
$\gln(K)$-$\gln(\hat\OO)$-orbit of $r$, so the map is injective.

To prove surjectivity, take an $\OO$-lattice $L\subset K^n_\infty$.
We have $KL=\Q L\cong\Q^{n[K\colon\Q]}$, so $\dim_K KL=n$. In
particular, $L$ is a finitely generated $\OO$-module of rank $n$.
Therefore it suffices to show that there exists $g\in\gln(K_\infty)$
such that $gL\subset K^n$. Let $e_1,\dots,e_n$ be a basis of $KL$
over $K$. Since $KL=\Q L$ is dense in $K_\infty^n$, the image of
$KL$ under the projection $K^n_\infty\to K^n_v$ is dense in $K^n_v$
for any infinite place $v$. It follows that the images of
$e_1,\dots,e_n$ are linearly independent over $K_v$. So there exists
$g_v\in\gln(K_v)$ which maps these images onto the standard basis of
$K_v^n$. Then $g=(g_v)_{v|\infty}$ is an element in $\gln(K_\infty)$
mapping $e_1,\dots,e_n$ onto the standard basis of $K^n_\infty$, so
that $gKL=K^n$. \ep

For $s\in\gln(\ak)$ and $t\in\mn(\hat\OO)$ consider the
$\OO$-lattice $L=s_f\hat\OO^n\cap K^n$. The map
$s_ft\colon\akf^n\to\akf^n$ maps $\hat\OO^n$ into $s_f\hat\OO^n$,
hence induces an $\hat\OO$-module map $\akf^n/\hat\OO^n\to
\akf^n/s_f\hat\OO^n$. Then there exists a unique $\OO$-module map
$\varphi\colon K^n/\OO^n\to KL/L$ such that the diagram
$$
\xymatrix{\akf^n/\hat\OO^n\ar[r]^{s_ft} & \akf^n/s_f\hat\OO^n\\
K^n/\OO^n\ar[u] \ar[r]_{\varphi} & KL/L\ar[u]}
$$
commutes, where the vertical arrows are the canonical isomorphisms
defined by the inclusions $K^n\subset\akf^n$, $KL\subset\akf^n$. We
shall also denote $\varphi$ by $[s_ft]$. Thus $(L,\varphi)$ is a
$K$-lattice. Therefore $(s_\infty^{-1}L,s_\infty^{-1}\varphi)$ is
also a $K$-lattice, which we denote by $[(s,t)]$.

\begin{lemma} \label{Klattice}
The map $\gln(\ak)\times\mn(\hat\OO)\ni
(s,t)\to[(s,t)]=(s_\infty^{-1}(s_f\hat\OO^n\cap
K^n),s_\infty^{-1}[s_ft])$ induces a bijection between
$$
\gln(K)\backslash\gln(\ak)\times_{\gln(\hat\OO)}\mn(\hat\OO)
$$
and the set of $n$-dimensional $K$-lattices.
\end{lemma}

\bp By Lemma~\ref{olat} we only need to check that any $\OO$-module
map $\akf^n/\hat\OO^n\to\akf^n/s_f\hat\OO^n$, where
$s_f\in\gln(\akf)$, is defined by the matrix $s_ft$ for a unique
$t\in\mn(\hat\OO)$. It suffices to consider $s_f=1$. The problem
then reduces to showing that any $\OO$-module map $K_v/\OO_v\to
K_v/\OO_v$ is given by multiplication by a unique element of
$\OO_v$. If $\pi$ is a uniformizing parameter in $\OO_v$, then any
$\OO$-module map $\OO_v\pi^{-m}/\OO_v\to K_v/\OO_v$ is determined by
the image of $\pi^{-m}$, so it is given by multiplication by an
element in $\OO_v$ which is uniquely determined modulo $\OO_v\pi^m$.
Since $\OO_v$ is complete in the $(\pi)$-adic topology, this gives
the result. \ep

Notice that we have shown in particular that for any $K$-lattice
$(L,\varphi)$ with $L\subset K^n$ the homomorphism $\varphi$ lifts
to a unique $\akf$-module map $\tilde\varphi\colon\akf^n\to\akf^n$.

\begin{definition} \label{defComm}
Two $n$-dimensional $K$-lattices $(L_1,\varphi_1)$ and
$(L_2,\varphi_2)$ are called commensurable if the lattices $L_1$ and
$L_2$ are commensurable and $\varphi_1=\varphi_2$ modulo $L_1+L_2$.
\end{definition}

If $L_1$ and $L_2$ are commensurable then $KL_1=\Q L_1=\Q L_2=KL_2$.
In particular, if $L_1\subset K^n$ then also $L_2\subset K^n$. It is
clear that then the lifting of the composition of the homomorphisms
$\varphi_1\colon K^n/\OO^n\to KL_1/L_1$ and $KL_1/L_1\to
K(L_1+L_2)/(L_1+L_2)$ coincides with $\tilde\varphi_1$. Therefore
two $K$-lattices $(L_1,\varphi_1)$ and $(L_2,\varphi_2)$ with
$L_1,L_2\subset K^n$ are commensurable if and only if
$\tilde\varphi_1=\tilde\varphi_2$. This implies that
commensurability is an equivalence relation.

Denote the equivalence relation of commensurability of
$n$-dimensional $K$-lattices by $\rr_{K,n}$. Consider now the action
of $\gln(\akf)$ on $\gln(K)\backslash\gln(\ak)\times\mn(\akf)$
defined by
$$
g(s,t)=(sg^{-1},gt).
$$
Define a subgroupoid
$$
\gln(\akf)\boxtimes(\gln(K)\backslash\gln(\ak)\times\mn(\hat\OO))
=\{(g,s,t)\mid t\in\mn(\hat\OO),\ gt\in\mn(\hat\OO)\}
$$
of the transformation groupoid
$\gln(\akf)\times(\gln(K)\backslash\gln(\ak)\times\mn(\akf))$. We
have a groupoid homomorphism
$$
\gln(\akf)\boxtimes(\gln(K)\backslash\gln(\ak)\times\mn(\hat\OO))\to
\rr_{K,n}
$$
defined by
\begin{equation} \label{ehomo}
(g,s,t)\mapsto ([(sg^{-1},gt)],[(s,t)]).
\end{equation}
To see that $[(s,t)]$ and $[(sg^{-1},gt)]$ are indeed commensurable
recall that by definition we have
$[(s,t)]=(s_\infty^{-1}(s_f\hat\OO^n\cap K^n),s_\infty^{-1}[s_ft])$
and $[(sg^{-1},gt)]=(s_\infty^{-1}(s_fg^{-1}\hat\OO^n\cap
K^n),s_\infty^{-1}[s_ft])$.

By Lemma~\ref{Klattice} to make the above homomorphism injective we
have to factor out the action of $\gln(\hat\OO)$. Consider the
action of $\gln(\hat\OO)\times\gln(\hat\OO)$ on
$\gln(\akf)\boxtimes(\gln(K)\backslash\gln(\ak)\times\mn(\hat\OO))$
defined by
$$
(u_1,u_2)(g,s,t)=(u_1gu_2^{-1},su_2^{-1},u_2t),
$$
and denote by
$$
\gln(\hat\OO)\backslash\gln(\akf)\boxtimes_{\gln(\hat\OO)}
(\gln(K)\backslash\gln(\ak)\times\mn(\hat\OO))
$$
the quotient space.

\begin{proposition}
The map (\ref{ehomo}) induces a bijection between
$$
\gln(\hat\OO)\backslash\gln(\akf)\boxtimes_{\gln(\hat\OO)}
(\gln(K)\backslash\gln(\ak)\times\mn(\hat\OO))
$$
and $\rr_{K,n}$.
\end{proposition}

\bp By Lemma~\ref{Klattice} the map
$$
\gln(\hat\OO)\backslash\gln(\akf)\boxtimes_{\gln(\hat\OO)}
(\gln(K)\backslash\gln(\ak)\times\mn(\hat\OO))\to\rr_{K,n}
$$
is well-defined and injective. To prove surjectivity we have to show
that if $(L,\varphi)=[(s,t)]$ is a $K$-lattice then any
commensurable $K$-lattice is of the form $[(sg^{-1},gt)]$ for some
$g\in\gln(\akf)$. We may assume that $L\subset K^n$ and then that
$s_\infty=1$. Then by Lemma~\ref{Klattice} and the discussion
following Definition~\ref{defComm} any commensurable $K$-lattice is
of the form $[(q,r)]$ with $q_\infty=1$ and $q_fr=s_ft$. Letting
$g=q_f^{-1}s_f$ we get $(q,r)=(sg^{-1},gt)$. \ep

\begin{remark}
In the case $K=\Q$, or more generally for fields with class number
one, there is a better description due to the fact that any
$\Z$-lattice is free. Indeed, by freeness we have
$\gln(\af)=\glnq\gln(\Zhat)$, where $\glnq$ is the group of rational
matrices with positive determinant. It follows that any
$\gln(\Zhat)\times\gln(\Zhat)$-orbit in
$\gln(\af)\times(\gln(\A)\times\mn(\Zhat))$ has a representative in
$\glnq\times((\gln(\R)\times\glnq)\times\mn(\Zhat))$. Furthermore,
the map
$$
\glnr\times\glnq\to\glnr,\ \ (g,h)\mapsto h^{-1}g,
$$
induces a bijection between $\glnq\backslash(\glnr\times\glnq)$ onto
$\glnr$. One may then conclude that $\rr_{\Q,n}$ can be identified
with
$$
\slnz\backslash\glnq\boxtimes_{\slnz}(\glnr\times\mn(\Zhat)),
$$
where the action of $\slnz\times\slnz$ on
$\glnq\times\glnr\times\mn(\Zhat)$ is given by
$$
(\gamma_1,\gamma_2)(g,h,m)=(\gamma_1g\gamma_2^{-1},\gamma_2h,\gamma_2m).
$$
\end{remark}

Consider now the case $n=1$ (and $K$ arbitrary). Then we conclude
that there is a bijection between $\rr_{K,1}$ and the subgroupoid
$$
(\akf^*/\ohs)\boxtimes((\ak^*/K^*)\times_\ohs\hat\OO)
$$
of the transformation groupoid
$(\akf^*/\ohs)\times((\ak^*/K^*)\times_\ohs\akf)$. We have an
action, called the scaling action, of $K^*_\infty$ on $K$-lattices:
if $(L,\varphi)$ is a $K$-lattice and $k\in K_\infty^*$ then
$k(L,\varphi)=(kL,k\varphi)$. It defines an action of $K^*_\infty$
on $\rr_{K,1}$. In our transformation groupoid picture of
$\rr_{K,1}$ it corresponds to the action of $K^*_\infty$ by
multiplication on $\ak^*/K^*$. Denote by $(K^*_\infty)^\circ$ the
connected component of the identity in $K^*_\infty$. Then we get the
following result.

\begin{corollary} \label{Klatscal}
The quotient of the equivalence relation $\rr_{K,1}$ of
commensurability of $1$-dimensional $K$-lattices by the scaling
action of the connected component of the identity in $K^*_\infty$ is
a groupoid that is isomorphic to
$$
(\akf^*/\ohs)\boxtimes((\ak^*/K^*(K^*_\infty)^\circ)
\times_\ohs\hat\OO).
$$
\end{corollary}

Recalling that $\akf^*/\ohs\cong J_K$ and
$\ak^*/\overline{K^*(K^*_\infty)^\circ} \cong \gal(\kab/K)$ by class
field theory, we see that the above groupoid is almost the same that
we used to define the BC-system. The small nuance is that when we
put $\gal(\kab/K)$ in our topological groupoid in Section~\ref{sBC}
we were effectively taking  the quotient of $\ak^*$ by  the {\em
closure} of  $K^*(K_\infty^*)^\circ$. In terms of $K$-lattices this
means that given a $K$-lattice $(L,\varphi)$ we would have to
identify not only all $K$-lattices $(kL,k\varphi)$ with $k\in
(K^*_\infty)^\circ$, but also all $K$-lattices of the form
$(kL,k\psi)$, where $\psi$ is a limit point of the maps $u\varphi$
with $u\in\OO^*\cap(K_\infty^*)^\circ$ in the topology of pointwise
convergence.

\end{document}